\output={\if N\header\headline={\hfill}\fi
\plainoutput\global\let\header=Y}
\magnification\magstep1
\tolerance = 500
\hsize=14.4true cm
\vsize=22.5true cm
\parindent=6true mm\overfullrule=2pt
\newcount\kapnum \kapnum=0
\newcount\parnum \parnum=0
\newcount\procnum \procnum=0
\newcount\nicknum \nicknum=1
\font\ninett=cmtt9

\font\ninebf=cmbx9

\font\sixbf=cmbx6
\font\ninesl=cmsl9

\font\nineit=cmti9

\font\ninerm=cmr9

\font\sixrm=cmr6
\font\ninei=cmmi9
\font\eighti=cmmi8
\font\sixi=cmmi6
\skewchar\ninei='177 \skewchar\eighti='177 \skewchar\sixi='177
\font\ninesy=cmsy9
\font\eightsy=cmsy8
\font\sixsy=cmsy6
\skewchar\ninesy='60 \skewchar\eightsy='60 \skewchar\sixsy='60
\font\titelfont=cmr10 scaled 1440
\font\paragratit=cmbx10 scaled 1200

\font\name=cmcsc10
\font\emph=cmbxti10

\font\tenmsbm=msbm10
\font\sevenmsbm=msbm7
%

%

%
\font\teneufm=eufm10
\font\seveneufm=eufm7
\font\fiveeufm=eufm5
\newfam\eufmfam
\textfont\eufmfam=\teneufm
\scriptfont\eufmfam=\seveneufm
\scriptscriptfont\eufmfam=\fiveeufm

\def\sl{{\rm slope}}
\font\tenmsam=msam10
\font\sevenmsam=msam7
\font\fivemsam=msam5
\newfam\msamfam
\textfont\msamfam=\tenmsam
\scriptfont\msamfam=\sevenmsam
\scriptscriptfont\msamfam=\fivemsam
\font\tenmsbm=msbm10
\font\sevenmsbm=msbm7
\font\fivemsbm=msbm5
\newfam\msbmfam
\textfont\msbmfam=\tenmsbm
\scriptfont\msbmfam=\sevenmsbm
\scriptscriptfont\msbmfam=\fivemsbm
\def\Bbb#1{{\fam\msbmfam\relax#1}}
\def\cz{{\kern0.4pt\Bbb C\kern0.7pt}
}
\def\az{{\kern0.4pt\Bbb A\kern0.7pt}
}
\def\ez{{\kern0.4pt\Bbb E\kern0.7pt}
}
\def\fz{{\kern0.4pt\Bbb F\kern0.3pt}}
\def\gz{{\kern0.4pt\Bbb Z\kern0.7pt}}
\def\hz{{\kern0.4pt\Bbb H\kern0.7pt}
}
\def\kz{{\kern0.4pt\Bbb K\kern0.7pt}
}
\def\nz{{\kern0.4pt\Bbb N\kern0.7pt}
}
\def\oz{{\kern0.4pt\Bbb O\kern0.7pt}
}
\def\rz{{\kern0.4pt\Bbb R\kern0.7pt}
}
\def\sz{{\kern0.4pt\Bbb S\kern0.7pt}
}
\def\pz{{\kern0.4pt\Bbb P\kern0.7pt}
}
\def\qz{{\kern0.4pt\Bbb Q\kern0.7pt}
}
\newskip\ttglue
\def\ninepoint{\def\rm{\fam0\ninerm}%
  \textfont0=\ninerm \scriptfont0=\sixrm \scriptscriptfont0=\fiverm
  \textfont1=\ninei \scriptfont1=\sixi \scriptscriptfont1=\fivei
  \textfont2=\ninesy \scriptfont2=\sixsy \scriptscriptfont2=\fivesy
  \textfont3=\tenex \scriptfont3=\tenex \scriptscriptfont3=\tenex
  \def\it{\fam\itfam\nineit}%
  \textfont\itfam=\nineit
  \def\sl{\fam\slfam\ninesl}%
  \textfont\slfam=\ninesl
  \def\bf{\fam\bffam\ninebf}%
  \textfont\bffam=\ninebf \scriptfont\bffam=\sixbf
   \scriptscriptfont\bffam=\fivebf
  \def\tt{\fam\ttfam\ninett}%
  \textfont\ttfam=\ninett
  \tt \ttglue=.5em plus.25em minus.15em
  \normalbaselineskip=11pt
  \font\name=cmcsc9
  \let\sc=\sevenrm
  \let\big=\ninebig
  \setbox\strutbox=\hbox{\vrule height8pt depth3pt width0pt}%
  \normalbaselines\rm
  \def\sl{\it}}

\headline={\ifodd\pageno\rightheadline\else\leftheadline\fi}
\def\rightheadline{\ninepoint Paragraphen"uberschrift\hfill\folio}
\def\leftheadline{\ninepoint\folio\hfill Chapter"uberschrift}
\let\header=Y
\def\titel#1{\need 9cm \vskip 2truecm
\parnum=0\global\advance \kapnum by 1
{\baselineskip=16pt\lineskip=16pt\rightskip0pt
plus4em\spaceskip.3333em\xspaceskip.5em\pretolerance=10000\noindent
\titelfont Chapter \uppercase\expandafter{\romannumeral\kapnum}.
#1\vskip2true cm}\def\leftheadline{\ninepoint
\folio\hfill Chapter \uppercase\expandafter{\romannumeral\kapnum}.
#1}\let\header=N
}
\def\Titel#1{\need 9cm \vskip 2truecm
\global\advance \kapnum by 1
{\baselineskip=16pt\lineskip=16pt\rightskip0pt
plus4em\spaceskip.3333em\xspaceskip.5em\pretolerance=10000\noindent
\titelfont\uppercase\expandafter{\romannumeral\kapnum}.
#1\vskip2true cm}\def\leftheadline{\ninepoint
\folio\hfill\uppercase\expandafter{\romannumeral\kapnum}.
#1}\let\header=N
}
\def\need#1cm {\par\dimen0=\pagetotal\ifdim\dimen0<\vsize
\global\advance\dimen0by#1 true cm
\ifdim\dimen0>\vsize\vfil\eject\noindent\fi\fi}
\def\neupara#1{\par\penalty-2000
\procnum=0\global\advance\parnum by 1
\vskip1cm\noindent{\paragratit \the\parnum. #1}%
\def\rightheadline{\ninepoint\S\the\parnum.\ #1\hfill \folio}%
\vskip 8mm\noindent}
\def\Proclaim #1 #2\finishproclaim {\bigbreak\noindent
{\bf#1\unskip{}. }{\it#2}\medbreak\noindent}
%
\gdef\proclaim #1 #2 #3\finishproclaim {\bigbreak\noindent%
\global\advance\procnum by 1
{%
{\relax\ifodd \nicknum
\hbox to 0pt{\vrule depth 0pt height0pt width\hsize
   \quad \ninett#3\hss}\else {}\fi}%
\bf\the\parnum.\the\procnum\ #1\unskip{}. }
{\it#2}
\immediate\write\num{\string\def
 \expandafter\string\csname#3\endcsname
 {\the\parnum.\the\procnum}}
\medbreak\noindent}
\newcount\stunde \newcount\minute \newcount\hilfsvar
\def\uhrzeit{
    \stunde=\the\time \divide \stunde by 60
    \minute=\the\time
    \hilfsvar=\stunde \multiply \hilfsvar by 60
    \advance \minute by -\hilfsvar
    \ifnum\the\stunde<10
    \ifnum\the\minute<10
    0\the\stunde:0\the\minute~Uhr
    \else
    0\the\stunde:\the\minute~Uhr
    \fi
    \else
    \ifnum\the\minute<10
    \the\stunde:0\the\minute~Uhr
    \else
    \the\stunde:\the\minute~Uhr
    \fi
    \fi
    }

 \def\calH{{\cal H}}
 
\def\calK{{\cal K}}

\def\dim{\mathop{\rm dim}\nolimits}

\def\GL{\mathop{\rm GL}\nolimits}

\def\kernel{\mathop{\rm kernel}\nolimits}

\def\mod{\mathop{\rm mod}\nolimits}

\def\Sp{\mathop{\rm Sp}\nolimits}

\def\boxit#1{
  \vbox{\hrule\hbox{\vrule\kern6pt
  \vbox{\kern8pt#1\kern8pt}\kern6pt\vrule}\hrule}}
\def\Boxit#1{
  \vbox{\hrule\hbox{\vrule\kern2pt
  \vbox{\kern2pt#1\kern2pt}\kern2pt\vrule}\hrule}}

\def\zwischen#1{\bigbreak\noindent{\bf#1\medbreak\noindent}}

\def\bigni{\bigskip\noindent }

\def\lo{\longrightarrow}

\def\betr#1{\vert#1\vert}

\def\pii{\pi {\rm i}}

\def\set#1{\bigl\{\,#1\,\bigr\}}

\def\square{\hbox{\hbox to 0pt{$\sqcup$\hss}\hbox{$\sqcap$}}}
\def\qed{\ifmmode\square\else{\unskip\nobreak\hfil
\penalty50\hskip3em\null\nobreak\hfil\square
\parfillskip=0pt\finalhyphendemerits=0\endgraf}\fi}
\def\pn{\the\parnum.\the\procnum}
\def\downmapsto{{\buildrel
        {\vbox{\hbox{\hskip.2pt$\scriptstyle-$}}}
        \over{\raise7pt\vbox{\vskip-4pt\hbox{$\textstyle\downarrow$}}}}}
\input level4_8.num
\nopagenumbers

\def\tr{{\rm tr}}

\def\eps{{\varepsilon}}

\immediate\newwrite\num
\nicknum=0  

\let\header=N
\immediate\newwrite\num\immediate\openout\num=level4_8.num
\def\RAND#1{\vskip0pt\hbox to 0mm{\hss\vtop to 0pt{%
\raggedright\ninepoint\parindent=0pt%
\baselineskip=1pt\hsize=2cm #1\vss}}\noindent}
\noindent
\centerline{\titelfont Siegel modular forms of level (4,8) and weight two}%
\def\leftheadline{\ninepoint\folio\hfill
Siegel modular forms of level (4,8) and weight two}%
\def\rightheadline{\ninepoint Introduction\hfill \folio}%
\headline={\ifodd\pageno\rightheadline\else\leftheadline\fi}
\vskip 1.5cm
\leftline{\it \hbox to 6cm{Eberhard Freitag\hss}
Riccardo Salvati
Manni  }
  \leftline {\it  \hbox to 6cm{Mathematisches Institut\hss}
Dipartimento di Matematica, }
\leftline {\it  \hbox to 6cm{Im Neuenheimer Feld 288\hss}
Piazzale Aldo Moro, 2}
\leftline {\it  \hbox to 6cm{D69120 Heidelberg\hss}
 I-00185 Roma, Italy. }
\leftline {\tt \hbox to 6cm{freitag@mathi.uni-heidelberg.de\hss}
salvati@mat.uniroma1.it}
\vskip1cm
\centerline{Heidelberg-Roma 2025}
\vskip1cm\noindent
{\paragratit Introduction}%
\vskip1cm\noindent 
We consider the space $[\Gamma_g[4,8],2]$ of Siegel modular forms
of genus $g$ of weight two. Examples of this space are the products
of 4 classical  theta nullwerte
$\vartheta[m]$ where
\def\d{\displaystyle}
$$\vartheta[m](\tau)=
\d\sum_{p\in\gz^g}\exp\bigl(\pii\tau[p+a]+
2(p+a)'b\bigr),\quad m=
\pmatrix{a\cr b}\in {1\over 2}\gz^{2g}.$$
One of the main results of this paper is that in the case $g\ge 8$ the space
$[\Gamma_g[4,8],2]$ is generated by the products of 4 theta nullwerte. 
We will obtain this as an application of the theory of singular modular forms [Fr].
We expect that this method carries over to $g\ge 5$. In the cases $g=1,2$ this result 
is also known [Ig]. The cases $g=3,\dots ,7$ remain open.
\smallskip
From this and the results of [SM1] we will get
$$\dim [\Gamma_g[4,8],2]=  {2^{g-1}(2^g +1)+3\choose 4  }- 
\sum_{i=0}^2 \mu_i (\nu_i -\pi_i)  \quad\hbox{for}\ g\ne  3,4,5, 6, 7$$
with  
$$\eqalign{
\mu_0&=1,\cr
\mu_i&=  \prod_{\nu=1}^i(2^{2(g-\nu+1)}-1)/ ( 2^{i-\nu+1}-1),\quad 0<i\le 2,\cr
\nu_i&= 2^{g-i-1}(2^{g-i}+1),\quad 0\le i\le2,\cr
\pi_i&=  (2^{g-i}+1)(2^{g-i-1}+1)/3,\quad 0\le i\le2.\cr}$$
The theta nullwerte are special cases of theta series which are attached
to positive definite real matrices. So, let $S$ be a positive 
definite $r\times r$-matrix and let
$A,B$ be real $r\times g$-matrices, Then we can define
$$\vartheta^S\left[\matrix{A\cr B}\right](\tau)=
\sum_{G\;{\rm integral}}\exp\pii \tr\bigl(S[G+A]\tau+2(G+A)'B\bigr).$$
If we specialize this to $S=(1)$, we get
$$\vartheta\left[\matrix{a\cr b}\right]=\vartheta^{(1)}
\left[\matrix{a'\cr b'}\right].$$
There is an easy generalization of this formula.
Let $A,B$ be two $r\times g$ matrices with entries in $\gz/2$. Denote by
$a_i,b_i$ there columns. Then
$$\prod_{i=1}^r\vartheta\left[\matrix{a_i\cr b_i}\right]=\vartheta^{E_r}
\left[\matrix{A'\cr B'}\right].$$
This  note can be considered as a  completion of the example  at the  end  of [Fr].
\zwischen{Mumford's theta relation}%
Let $S,T$ be two rational positive definite $r\times r$-matrices and let $A$ be a 
rational matrix such that
$$S=T[A].$$
Consider the finite groups
$$\eqalign{
\calK_1&=A\gz^{r\times g}/(A\gz^{r\times g}\cap\gz^{r\times g}),\cr
\calK_2&=A'^{-1}\gz^{r\times g}/(A'^{-1}\gz^{r\times g}\cap\gz^{r\times g}).\cr}$$
Then for any two rational $r\times g$-matrices $P,Q$ the relation
$$\vartheta^S\left[\matrix{A^{-1}P\cr A'Q}\right]={1\over\#\calK_2}
\sum_{X\in\calK_1,\;Y\in\calK_2}e^{-2\pii\tr(P'Y)}
\vartheta^T\left[\matrix{P+X\cr Q+Y}\right]$$
holds [Mu], Theorem 6.1.
\smallskip
Besides the theory of singular modular forms, this relation will play an important
role for the proof.
 
\neupara{Notations and Definitions}%
We use the following notations for matrices $A$. 
The transposed of $A$ is denoted
by $A'$.  Its trace is denoted by $\tr(A)$ and $A_0$ is the 
diagonal
of $A$ written
as column vector.
 Let $A$ be an $n\times n$-matrix and $B$ an $n\times m$-matrix, then
$$A[B]:=B'AB.$$
We have to consider certain congruence subgroups of the Siegel modular group
$$\Gamma_g=\Sp(g,\gz)\subset\GL(2g,\gz),$$
namely the principal congruence subgroup
$$\Gamma_g[q]=\kernel(\Sp(g,\gz)\lo \Sp(g,\gz/q\gz))$$
and Igusa's congruence group
$$\Gamma_g[q,2q]=\set{M\in \Gamma_g[q];\quad (AB')_0\equiv(CD')_0
\equiv0\mod 2q}.$$
Here $M={A\,B\choose C\,D}$ is the decomposition of $M$ into 4 blocks.
\smallskip
The Siegel half plane $\calH_g$ of genus $g$ is the set of all complex
$g\times g$-matrices such that its imaginary part is positive definite. The
modular group acts on $\calH_g$,
$$M\tau=(A\tau+B)(C\tau+D)^{-1}.$$
We choose for each $M\in\Gamma_g$ a holomorphic square root
$\sqrt{\det(C\tau+D)}$. 
Let $\Gamma\subset \Gamma_g$ be some congruence subgroup. 
A multiplier system on
$\Gamma$ of weight $r/2$, $r\in \gz$, is a function
$$v:\Gamma\lo S^1=\{\zeta;\ \betr\zeta=1\}$$
such that $v(M)\sqrt{\det(C\tau+D)}\,\strut^r$ is an automorphy factor. If $r$ is even, this
means that $v$ is a character.
\smallskip
A modular form on $\Gamma$ of weight $r/2$ and with respect to
the multiplier system $v$ is a holomorphic function
$$f:\Gamma\lo\cz,\quad f(M\tau)=v(M)
\sqrt{\det(C\tau+D)}\,\strut^r\ (M\in\Gamma),$$
where in the case $g=1$ the usual regularity condition at the cusps has to
be added. The space of these forms is denoted by
$$[\Gamma,r/2,v]$$
If $r$ is even and $v$ is trivial we omit $v$ in the notation.
\neupara{Isotropic matrices and theta series}%
We  consider the group $\Gamma_g[q,2q]$. We have to  consider 
positive definite
$r\times r$-matrices $S$ such that $S$ and $q^2S^{-1}$ are integral. A 
$g\times g$ matrix $V$ is called isotropic for $(S,q)$ if
$$S^{-1}[V]\quad\hbox{and}\quad qS^{-1}V$$
are integral. We identify elements of $(\gz^r)^g$ with $r\times g$ matrices. The 
condition that $V$ is isotropic is a condition mod
$$(q\gz^r+S\gz^r)^g.$$
Hence we can consider isotropic matrices as elements of
$$(\gz^r/(q\gz^r+S\gz^r))^g.$$
One knows ([Fr], Corollary II.6.11) that the theta 
series
$$\vartheta_{S,V}(\tau)=\sum_{G\;{\rm integral}}
\exp{\pii\over q}\tr( S[G]\tau+2G'V)$$
is a modular form on  $\Gamma_g[q,2q]$ of weight $r/2$ and a certain multiplier
system $\eps_S$ which is independent of $V$
$$\vartheta_{S,V}\in[\Gamma_g[4,8],r/2,\eps_S].$$
In the notation of the intoduction we have
$$\vartheta_{S,V}=\vartheta^{S/q}\left[\matrix{0\cr V/q}\right].$$
\proclaim
{Definition}
{Let S be a positive definite $r\times r$-matrix such that $S$ and $q^2S^{-1}$ are
integral. The space $\Theta(S,g,q)$ is the span of all theta series
$\vartheta_{S,V}$ with isotropic $V$.}
DTh%
\finishproclaim
So we have
$$\Theta(S,g,q)\subset [\Gamma_g[4,8],r/2,\eps_S].$$
A group $L\subset \gz^r/(q\gz^r+S\gz^r)$ 
is called isotropic (with respect to $(S,q)$) if there exists an isotropic
$V$ such that $L$ is generated by the columns of $V$. 
A subgroup of an isotropic group is isotropic. 
\proclaim
{Definition}
{Let S be a positive definite $r\times r$-matrix such that $S$ and $q^2S^{-1}$ are
integral. Let $L\subset \gz^r/(q\gz^r+S\gz^r)$  be an isotropic subspace
The space $\Theta_L(S,g,q)$ is the span of all theta series
$\vartheta_{S,V}$ such that the columns of $V$ are contained in $L$.}
DThz%
\finishproclaim
So we have
$$\Theta(S,g,q)=\sum_{L\;{\rm isotropic}}\Theta_L(S,g,q).$$
One result of the theory singular modular forms is the following theorem.
\proclaim
{Theorem}
{Let $\eps$ be a multiplier system on $\Gamma_g[q,2q]$ such that there exists
a positive definite matrix $S$ such that $S$ and $q^2S^{-1}$ are integral and such that
$\eps=\eps_S$. Assume $g\ge 2r$. Then
$$ [\Gamma_g[4,8],r/2,\eps]=\sum_{\eps=\eps_S}\Theta(S,g,q).$$
Even more is true. It is sufficient to restrict in this sum to $S$ such that
$qS^{-1}$ is integral.}
Tq%
\finishproclaim
For the proof we refer to [Fr], Theorem VI.1.5 and Theorem VI.1.6. 
(See also Proposition \PTt. It gives a simple proof that Theorem VI.1.5 in [Fr]
implies Theorem VI.1.6.)
\qed\smallskip
In the general (not necessarily isotropic case) we define
$$ [\Gamma_g[4,8],r/2,\eps]_\Theta=\sum_{\eps=\eps_S}\Theta(S,g,q)$$
where $S$ runs through all $S$ as in the theorem.
\neupara{Multiplier systems}%
 \proclaim
 {Lemma}
 {Let $S$ be a positive definite $4\times 4$-matrix such that
 $S$ and $16S^{-1}$ are integral. The multiplier system $\eps_S$
 is trivial on $\Gamma_g[4,8]$ if and only if  the determinant of $S$ is a square.
}
Lisq%
\finishproclaim
{\it Proof.} 
Assume that $S$ and $16S^{-1}$ are integral. The determinant of $S$ is a power of 
$2$. We must show that it is an ever power of 2.
The multiplier system $\eps_S$  can be computed as follows. We can assume 
$V=0$,
$$\vartheta_{S,0}(\tau)=\sum_{G\;{{\rm integral}}}\exp {\pii\over 4} \tr (S[G]\tau)$$
Consider first
$$\vartheta_{4S,0}(\tau)=\sum_{G\;{{\rm integral}}}\exp {\pii} \tr (S[G]\tau)$$
It is known (i.e. [Fr], Proposition 7.1) that this is a modular form on the group
$$\Gamma_{g,0,\vartheta}[16]=\set{M\in\Gamma_g;\quad C\equiv 0\mod 16,\quad
\hbox{the diagonal of}\  (CD')/16\ \hbox{is even}}.$$
The multiplier system on this group is
$$\left({\det S\over\vert\det D\vert}\right)\qquad
\hbox{(generalized Legendre  symbol)}.$$
Now use
$$\vartheta_{S,0}(\tau)=\vartheta_{4S,0}(\tau/4).$$
This implies
$$\vartheta_{S,0}(M\tau)=\vartheta_{4S,0}(M(\tau)/4).$$
We have 
$$\eqalign{
M(\tau)/4=&\pmatrix{E/2&0\cr0&2E}M(\tau)\cr
=&\pmatrix{E/2&0\cr0&2E}M\pmatrix{2E&0\cr0&E/2}(\tau/4)\cr
=&\pmatrix{A&B/4\cr 4C&D}(\tau/4)\cr}$$
Assume $M\in \Gamma_g[4,8]$. Then
$$N=\pmatrix{A&B/4\cr 4C& D}\in \Gamma_{g,0,\vartheta}[16].$$
 We obtain 
 $$\vartheta_{4S,0}(N\tau)=
 \left({\det S\over\vert\det D\vert}\right)\det(C\tau+D)^2
 \vartheta_{4S,0}(\tau/4)$$
 or
 $$\vartheta_{S,0}(M\tau)=
 \left({\det S\over\vert\det D\vert}\right)\det(C\tau+D)^2
 \vartheta_{S,0}(\tau).$$
 This means
 $$\eps_S(M)= \left({\det S\over\vert\det D\vert}\right).$$
 We choose an element of $\Gamma_g[4,8]$ such that $\det D=5$. In the case
$g=1$ one can take  ${13\;8\choose 8\;\;5}$. This element can be embedded into
$\Gamma_g[4,8]$. Since $\bigl({2\over 5}\bigr)=-1$, the determinant of 
$S$ must be an  even power of $2$.
\qed
\neupara{Mumford's theta relation}%
In this section we treat three  examples for Mumford's theta relations.
\smallskip
The theta nullwerte
$\vartheta[m]$
are modular forms on $\Gamma_g[4,8]$ of weight $1/2$ with respect to a
joint multiplier system $v_\vartheta$. The square of it is trivial. Hence
the products of $r$ theta nullwerte are contained in 
$[\Gamma_g[4,8],r/2,v_\vartheta^r]$.
Actually
$$v_\vartheta^r=v_S\ \hbox{where}\ S=E_r\  ( r\times r\hbox{-unit matrix}).$$
The following proposition is the first example of Mumford's theta relations.
\proclaim
{Proposition}
{The space $\Theta(E_r,g,4)$ contains the monomials of $r$ theta null\-we\-rte
and is generated by them.}
Ptg%
\finishproclaim
{\it Proof.} Consider
$$\vartheta_{E,V}\in \Theta(E_r,g,4),$$
i.e.~
$$\vartheta_{E,V}(\tau)=\sum_G\exp\pi/4(\tr(E[G]\tau+2V'G)$$
with isotropic $V$. Isotropy in this case means simply that $V$ is integral.
Obviously 
$$\vartheta_{E,V}=\vartheta^{E/4}\left[\matrix{0\cr V/4}\right].$$
We apply Mumford's relation for $S=E/4$, $T=E$ and $A=E/2$.
$$\vartheta^{E/4}\left[\matrix{0\cr V/4}\right]=
\sum_{X\in\calK_1,\;Y\in\calK_2}e^{-\pii\tr(X'Y)}
\vartheta^E\left[\matrix{X\cr V/2+Y}\right]$$
Since $2X$ and $V+2Y$ are integral, the functions
$$\vartheta^E\left[\matrix{X\cr V/2+Y}\right]$$
are products of theta constants. 
\smallskip
Viceversa the trivial  fact that   $E =2 (E/4) 2$, we can write
$$\vartheta^E\left[\matrix{A\cr B}\right]=  \vartheta^{2 (E/4) 2}    
\left[\matrix{\alpha /2\cr 2\beta}\right] $$
with $\alpha$ integral  and $4\beta$ integral, thus, applying  Mumford's 
formula,  with $A=2E$,
we will sum over
$$\eqalign{
D\in\calK_1&=2\gz^{r\times g}/(2\gz^{r\times g}\cap\gz^{r\times g}),\cr
C\in \calK_2&=2^{-1}\gz^{r\times g}/(2^{-1}\gz^{r\times g}
\cap\gz^{r\times g}).\cr}$$
This  means $D=0$, $C$ half integral and   $d= 2^{gh}$. Hence we have
$$\vartheta^E\left[\matrix{A\cr B}\right]= 
2^{-gh} \sum_{ C\,half \, integral }  e(- tr(C'\alpha))  
\vartheta^{E/4} \left[\matrix{\alpha\cr  \beta  +C}\right].$$
We observe that  since  the  matrix  $\alpha$  is integral,
$$  \vartheta^{E/4} \left[\matrix{\alpha\cr  \beta  +C}\right]=  
\vartheta^{E/4} \left[\matrix{0\cr  \beta  +C}\right]$$
Moreover  the   matrices   $4(\beta+C)$ are  integral, hence 
they are    isotropic with respect to $E$.
Thus  theta  constants are linear combinations of theta 
series  in $\Theta(E_r,g,4)$ and they span the space.
\qed \smallskip
With the same method one can show that for more $(S,V)$ the theta series 
$\vartheta_{S,V}$ can be expressed by theta nullwerte.
\smallskip
Here is a second example for Mumford's thete relations.
\proclaim
{Proposition}
{Let $S$ be a positive definite $r\times r$-matrix sich that
$S$ and $16S^{-1}$ are integral
Assume that for each 
isotropic
$V$ there exists a solution $S=A'A$  where $A$ is an $r\times r$-matrix with
the properties
$$A\ \hbox{integral},\quad 4A^{-1}\ \hbox{integral},\quad A'^{-1}V\ \hbox{integral}.$$
Then the space $\Theta(S,g,4)$ is contained in the space generated by monomials
of degree $r$ of the theta nullwerte.}
PSq%
\finishproclaim
{\it Proof.} 
Let $\vartheta_{S,V}\in \Theta(S,g,q)$
Recall
$$\vartheta_{S,V}=\vartheta^{S/4}\left[\matrix{0\cr V/4}\right].$$
Obviously $S/4= (A'/2)(A/2)$ with  $A/2$ and   $(A/2)^{-1}$   half integral. 
Thus in 
 $$\vartheta^{S/4}\left[\matrix{0\cr V/4}\right]=
 \vartheta^{S/4}\left[\matrix{0\cr  (A'/2)(A'/2)^{-1}V/4}\right]$$
 $(A'/2)^{-1}V/4$  is half integral  and  also  all  characteristics  in 
 $$\eqalign{
\calK_1&=(A/2)\gz^{r\times g}/(A/2)\gz^{r\times g}\cap\gz^{r\times g}),\cr
\calK_2&=2A'^{-1}\gz^{r\times g}/2A'^{-1}\gz^{r\times g}\cap
\gz^{r\times g}).\cr}$$
 are half integral, thus $\vartheta_{S,V}$
 is a linear combinations of monomials in the   theta nullwerte.
 \qed\smallskip
Now we treat the third example for Mumford's theta relations.
We consider again a positive $r\times r$-matrix such
that $S$ and $q^2S^{-1}$ are integral, Let
$$L\subset \gz^r/(q\gz^r+S\gz^r)$$ 
be an isotropic subgroup. 
We need also the natural map
$$\gz^r\lo ((\gz^r/(q\gz^r+S\gz)^r))/L.$$
Its kernel
is of the form $A\gz^r$ where $A$ is an integral $r\times r$ matrix. 
From $A\gz^r\supset q\gz^r+S\gz^r$ follows that besides $A$ also $qA^{-1}$ and
$A^{-1}S$ are integral.  The columns of $A$ considered mod  $q\gz+S\gz^r$
are contained in $L$. Hence
$$\tilde S=S^{-1}[A].$$
is integral.
The matrices $\tilde S$ and $q\tilde S^{-1}=(A^{-1}S)(qA^{-1})$ are integral. Hence 
$q^2\tilde S^{-1}$ is integral too. 
The matrix $H=A^{-1}SG$ is integral for integral $G$. 
\proclaim
{Proposition}
{With the notations above we have
$$\Theta_L(S,g,q)\subset\Theta_{\tilde L}(\tilde S,g,q).$$
}
PTt%
\finishproclaim
{\it Proof.}
We apply to 
$S/4=(\tilde S/4)[A^{-1}S]$ 
Mumford's theta relation.
$$ \vartheta^{S/q}\left[\matrix{0\cr V/q}\right]=  
\vartheta^{S/q}\left[\matrix{0\cr SA^{'-1}A' S^{-1}V/q}\right]=
\sum_{ C\in \calK_2} \vartheta^{\tilde S/q}\left[\matrix{0\cr 
 A' S^{-1}V/q  +C}\right]$$
$$
\calK_2=A' S^{-1}\gz^{r\times g}/(A' S^{-1}\gz^{r\times g}\cap\gz^{r\times g}).$$
 Now since $q A' S^{-1}=q\tilde SA^{-1}$ is integral, $qC= qA' S^{-1} G$ 
 is integral for any integral $G$ and  
 $$\tilde S^{-1}[A'S^{-1}V+qC]= \tilde S^{-1}[A'S^{-1}(V+q G)]= $$
 $$  S^{-1}A\tilde S^{-1}A'S^{-1}[(V+q G)]  =   
  S^{-1}[(V+q G)] $$  
 that is  integral , hence $A'S^{-1}V+qC$ is  isotropic. Thus
 $$\Theta_L(S,g,q)\subset \Theta_{\tilde L}(\tilde S,g,q).\eqno\square$$
\neupara{Applications}%
We are intersted in the space $[\Gamma_g[4,8],2]$ of modular forms of
weight 2 on the group $\Gamma_g[4,8]$. Recall that 
$$[\Gamma_g[4,8],2]_\Theta\subset [\Gamma_g[4,8],2]$$ 
is the subspace generated by theta series of the form $\vartheta_{S,V}$.
We know that this space contains the products of $4$ theta nullwerte.
\proclaim
{Theorem}
{The space $[\Gamma_g[4,8],2]_\Theta$ equals the space generated by
products of $4$ theta nullwerte. Its dimension equals
$$\dim [\Gamma_g[4,8],2]_\Theta=
 {2^{g-1}(2^g +1)+3\choose 4  }- 
\sum_{i=0}^2 \mu_i (\nu_i -\pi_i) ,$$
where $\mu_i,\nu_i,\pi_i$ are defined in the introduction.}
TTv%
\finishproclaim
The proof is  a consequence of the following lemma.
\proclaim
{Lemma}
{Let $S$ be positive definite $4\times 4$-matrix such that $S$ and $4S^{-1}$ are
integral and that the determinant of $S$ is a square. Then each isotropic matrix
$V$ for $(S,4)$ there exists an integral $4\times 4$-matrix  $A$ such that
$S=A'A$ and such that  $4A^{-1}$ and $A^{-1}$ are integral.
}
LStri%
\finishproclaim
{\it Proof.} The proof rests on computer calculations. In [Fr] one finds a list of
representatives of the unimodular classes of all $S$ such that 
$S$ and $16S^{-1}$ are
integral and that the determinant of $S$ is a square. This list contains 138
elements. If one singles out those ones that already $4S^{-1}$ is integral, one gets
16 matrices $S$. For each $S$ one can compute the maximal isotropic
aubspaces $L\subset (\gz/4\gz)^4$. Let $L$ be one of them. It can be described
by 4 generators. We write them as columns $v_1,\dots,v_4$
and consider the matrix $V=(v_1,\dots.v_4)$. Then each isotropic
matrix related to $L$ is of the form $VG$  with an integral matrix $G$. Next one
computes the set of all integral solutions $S=A'A$. One shows that this set is 
not empty and then one checks thet for each $V$ there exitsts an $A$ such that
$A$, $4A'^{-1}$, $A^{-1}$ are integral.
Instead of replicating the program here, we explain it for an example.
$$S=\pmatrix{2&0&1&1\cr0&2&1&-1\cr1&1&2&0\cr1&-1&0&2\cr}.$$
There are 3 maximal isotropic groups. Their defining matrix is
$$
V_1=\pmatrix{
0&0&0&1\cr
0&0&1&0\cr
0&1&0&0\cr
2&1&0&0\cr},
\quad
V_2=
\pmatrix{
0&0&0&1\cr
0&0&1&0\cr
0&1&0&0\cr
2&1&1&1\cr},
\quad
V_3=
\pmatrix{
0&0&0&1\cr
0&0&2&1\cr
0&1&0&0\cr
1&0&0&0\cr}.
$$
In terms of isotropic subgroups this means
$$\eqalign{
L_1&=\set{x\in (\gz/4\gz)^4;\quad x_3+x_4\equiv 2\mod 4},\cr
L_2&=\set{x\in (\gz/4\gz)^4;\quad x_1+x_2+x_3=x_4\equiv 2\mod 4},\cr
L_3&=\set{x\in (\gz/4\gz)^4;\quad  x_1+x_2\equiv 2\mod 4}.\cr
}$$
One computes three solutions $A_1,A_2,A_3 $ of the equation $S=A'A$ such that
$A_i,4A_i^{-1},A_i'^{-1}V_i$ are integral, namely
$$
A_1=\pmatrix{
0&1&0&0\cr
0&1&1&-1\cr
1&0&0&0\cr
1&0&1&1\cr},
\
A_2=\pmatrix{
0&1&0&-1\cr
0&1&1&0\cr
1&0&0&1\cr
1&0&1&0\cr},
\
A_3=\pmatrix{
0&0&0&1\cr
0&0&1&0\cr
1&1&1&0\cr
1&-1&0&1\cr}.
$$
\proclaim
{Theorem}
{Assume $g\ge 8$. Then $[\Gamma_g[4,8],2]$ equals the space generated
by the fourth products  theta  nullwerte.}
TH%
\finishproclaim
\neupara{Related cases}%
There are more cases that can be treated with the same method. We keep short.
\smallskip
We need the theta nullwerte of second kind [Ru1]
$$f_a(\tau)=\vartheta\left[\matrix{a/2\cr 0}\right](2\tau),\quad a\in\gz^g..$$
One knows that the the $f_af_b$ generate the same vector space as  the squares
of the theta nullwerte.
\proclaim
{Theorem}
{Assume $g\ge 6$. 
The products of 4 theta nullwerte of second kind $f_a$ give a basis of
$[\Gamma_g[2,4],2]$.
Hence the dimension of this space is 
$$\dim[\Gamma_g[2,4],2]=\pmatrix{2^g+3\cr 4}.$$}
Tz%
\finishproclaim
{\it Proof.} In the case $g\ge 8$ this is a consequence of  Theorem \TTv. 
This follows from the results in [SM1]. The dimension formula is in
[SM1], Theorem 1, (ii) in a slightly different form. It says that the products
of four $f_a$ are linearly independent.
\smallskip
Another way to prove it is to apply again [Fr]. One has to determine representatives
of the unimodular classes of positive definite $4\times4$-matrices such that
$S$ and $2S^-1$ are integral. One shows that there are 6 classes. 
\smallskip
This second proof shows that the Theorem is 
true for $g\ge 6$. In [Fr] it  is shown 
that the bound $g\ge 8$ sometimes can be replaced by $g\ge 6$. This is true
if for each positive integral $S$ such that $S$ and $qS^{-1}$ are integral the
set of isotropic matrices is a group.(see [Fr] Proposition V.3.1). One can check
this for the 6 representatives.\qed
\proclaim
{Theorem}
{Assume $g\ge 6$. The space $[\Gamma_g[2],2]$ is generated 
by the fourth powers 
$$\vartheta\left[\matrix{a\cr b}\right]^4,\quad a,b\in \{0,1/2\}^g.$$
Its dimension is
$$(2^g+1)(2^{g-1}+1)/3.$$}
Tvp%
\finishproclaim
{\it Proof.} That  $[\Gamma_g[2],2]$ is generated by the fourth powers of the
theta nullwerte can be derived from Theorem \Tz. The dimension formula
has been proved in [vG].
\qed\smallskip
In the cases $g=1,2,3$,  the structure of the  rings of modular forms for the groups
$[\Gamma_g[2]$ and $[\Gamma_g[2,4]$
are known, cf. [Ig] and [Ru1],  [Ru2]. Thus   we can  say also that    
Theorem \Tz\  and  Theorem \Tvp\
hold for $g=1,2,3$.  We conclude  with the  following result:
 \proclaim
{Theorem}
{For arbitrary $g$ there is no non-vanishing cuspform in $[\Gamma_g[2],2]$}
 Tsl%
\finishproclaim
{\it Proof.} in the cases $g\le 3$ this follows from the structure theorems for the
rings of modular forms ([Ig], [Ru1], [Ru2]). In the cases $g\ge 5$ the modular forms
of weight 2 are singular, hence never cusp forms. The case $g=4$
needs an extra argument.
\smallskip
Let 
$$f(\tau)=\sum_T a(T)\exp2\pii\tr(T\tau)$$
be the Fourier expansion of a non-vanishing
Siegel modular form of genus $g$ and weight
$k$ on some congruence group.
Here $T$ runs through the semipositive matrices of a lattice of rational symmetric
matrices. We define
$$v_\infty(f)=\min\set{t_{11};\quad a(T)\ne 0}.$$
For $M\in\Gamma_g$ we define the transformed form by
$$(f\vert M)(\tau)=\det(C\tau+D)^{-k}f(\tau).$$
We set
$$v(f)=\min_{M\in\Gamma_g} v_\infty(f\vert M).$$
The {\it slope\/} of $f$ is defined by
$$\sl(f)={k\over v(f)}\le \infty.$$
It is finite if $f$ is a cusp form. Now, let $f\in[\Gamma_g[2],2]$ be a non-vanishing
cusp form. Then $v_\infty(f)\ge 1/2$. Since $\Gamma_g[2]$ is 
normal in $\Gamma_g$, it follows $v(f)\ge 1/2$. Hence $\sl(f)\le 4$. In
[SM2] it has been proved that, in genus $g=4$,  $\sl(f)\ge 8$. Hence such a cusp form cannot exist.
\vskip1.5cm\noindent
{\paragratit References}%
\bigni
\item{[Fr]} Freitag, E.:
{\it Singular modular forms and theta relations,}
Lecture Notes {\bf 1487},
Spinger-Verlag, Berlin, Heidelberg u.a., 1991
\medskip
\item{[Ig]} Igusa, J.I.:
{\it On Siegel Modular Forms of Genus Two (II),}
Am. J. of Math. Vol. {\bf86}, No 2, 1964
\medskip
\item{[Mu]} Mumford, D.:
{\it Tata Lectures on Theta I,}
Modern Birkh\"auser Classics, Birkh\"auser Boston, 1990
\medskip
 \item{[Ru1]} Runge, B.:
 {\it On Siegel modular forms Part I,}
J. Reine Angew. Math. {\bf436}, 57--85,  1993.
\medskip \item{[Ru2]} Runge, B.:
 {\it On Siegel modular forms Part  II,}
Nagoya Math. J.  Vol {\bf 138}, 179--197,  1995.
\medskip
\item{[SM1]}  Salvati-Manni, R.:
{\it On the dimension of the vector space $\cz[\theta_m]_4$},
Nagoya Math. J. {\bf 98}, 99--107,  1985.
\medskip
\item{[SM2]}  Salvati Manni, R.:
{\it Modular Forms of the fourth
degree,}
LNM, Proceedings Trento {\bf 1515}, 106--111,  1992
\medskip
\item{[vG]} van Geemen, B.:
{\it Siegel modular forms vanishing on the moduli space of curves,}
Invent. math. {\bf 78}, 1984
\bye